\documentclass[11pt,reqno]{amsart}
\usepackage{amsmath,amssymb,amsthm}
\usepackage{graphicx}
\usepackage[all]{xypic}
\usepackage[utf8,nocaptions]{vietnam}
\usepackage{calrsfs}
\input xypic
\def\thmsection{section}
\def\thmchangesection{changesection}

\def\thmchangechapter{changechapter}
\def\thmchange{change}
\def\thmplain{plain}
\ifx\theoremnumberstyle\thmplain
  \theoremstyle{break-italic}
  \newtheorem{satz}{Satz}
\else
  \ifx\theoremnumberstyle\thmsection
    \theoremstyle{break-italic}
    \newtheorem{satz}{Satz}[section]
  \else
    \ifx\theoremnumberstyle\thmchange
      \swapnumbers
      \theoremstyle{break-italic}
      \newtheorem{satz}{Satz}
    \else
      \ifx\theoremnumberstyle\thmchangesection
         \swapnumbers
         \theoremstyle{break-italic}
         \newtheorem{satz}{Satz}[section]
      \else
        \ifx\theoremnumberstyle\thmchangechapter
           \swapnumbers
           \theoremstyle{break-italic}
           \newtheorem{satz}{Satz}[chapter]
        \else
          \ifx\theoremnumberstyle\thmsection
             \theoremstyle{break-italic}
             \newtheorem{satz}{Satz}[section]
          \else
            \theoremstyle{break-italic}
            \newtheorem{satz}{Satz}[section]
          \fi
        \fi
      \fi
    \fi
  \fi
\fi

\theoremstyle{break-italic}
\newtheorem{theorem}[satz]{Theorem}
\newtheorem{lemma}[satz]{Lemma}

\newtheorem{corollary}[satz]{Corollary}
\newtheorem{Proposition}[satz]{Proposition}

\newtheorem*{conjecture*}{Conjecture}

\theoremstyle{break-roman}
\newtheorem{definition}[satz]{Definition}

\newtheorem{example}[satz]{Example}

\newtheorem{remark}[satz]{Remark}

\newtheorem{conjecture}[satz]{Conjecture}

\theoremstyle{standard}

\newtheorem*{claim}{Claim}

\theoremstyle{varthm-roman}
\newtheorem*{varthm-roman}{}% all text supplied in the note
\theoremstyle{varthm-italic}
\newtheorem*{varthm-italic}{}% all text supplied in the note
\theoremstyle{varthm-roman-break}
\newtheorem*{varthm-roman-break}{}% all text supplied in the note
\theoremstyle{varthm-italic-break}
\newtheorem*{varthm-italic-break}{}% all text supplied in the note
\theoremstyle{varthm-roman-no-punctuation}
\newtheorem{varthm-roman-no-punctuation-numbered}[satz]{}% all text supplied in the note
\theoremstyle{varthm-italic-no-punctuation}
\newtheorem{varthm-italic-no-punctuation-numbered}[satz]{}% all text supplied in the note

\newenvironment{varthm-roman-numbered}[1]{
  \begin{varthm-roman-no-punctuation-numbered}
    \mbox{\rm\textbf{#1}}
  }{\end{varthm-roman-no-punctuation-numbered}}
\newenvironment{varthm-italic-numbered}[1]{
  \begin{varthm-italic-no-punctuation-numbered}
    \mbox{\rm\textbf{#1}}
  }{\end{varthm-italic-no-punctuation-numbered}}
\newenvironment{varthm-roman-break-numbered}[1]{
  \begin{varthm-roman-no-punctuation-numbered}
    \mbox{\rm\textbf{#1}\newline}
  }{\end{varthm-roman-no-punctuation-numbered}}
\newenvironment{varthm-italic-break-numbered}[1]{
  \begin{varthm-italic-no-punctuation-numbered}
    \mbox{\rm\textbf{#1}}\newline
  }{\end{varthm-italic-no-punctuation-numbered}}
\numberwithin{equation}{section}
\def\ex{\begin{example}
  }
  \def\eex{\end{example}}
\def\thr{\begin{theorem}}
\def\ethr{\end{theorem}}
\def\pro{\begin{Proposition}}
\def\epro{\end{Proposition}}
\def\coro{\begin{corollary}}
\def\ecoro{\end{corollary}}
\def\df{\begin{definition}}
\def\edf{\end{definition}}
\def\lm{\begin{lemma}}
\def\elm{\end{lemma}}
\def\pf{\begin{proof}}
\def\epf{\end{proof}}
\def\problem{\begin{problem}}
\def\eproblem{\end{problem}}

\def\ite{\begin{itemize}}
\def\hite{\end{itemize}}
\def\rem{\begin{remark}}
\def\erem{\end{remark}}

\def\cla{\begin{claim}}
\def\ecla{\end{claim}}
\def\conj{\begin{conjecture}}
\def\econj{\end{conjecture}}
\def\eex{{\accent"5E e}\kern-.385em\raise.2ex\hbox{\char'23}\kern-.08em}
\def\EES{{\accent"5E E}\kern-.5em\raise.8ex\hbox{\char'23 }}
\def\ow{o\kern-.42em\raise.82ex\hbox{
\vrule width .12em height .0ex depth .075ex \kern-0.16em \char'56}\kern-.07em}
\def\OW{O\kern-.460em\raise1.36ex\hbox{
\vrule width .13em height .0ex depth .075ex \kern-0.16em \char'56}\kern-.07em}
\begin{document}

\title[On the location of  eigenvalues of matrix polynomials]{On the location of   eigenvalues of matrix polynomials}

\author{L\^{e} C\^{o}ng-Tr\`{i}nh }
\address{L\^{e} C\^{o}ng-Tr\`{i}nh \\
Department of Mathematics, Quy Nhon University\\
170 An Duong Vuong, Quy Nhon, Binh Dinh}
%    Current address 
%\curraddr{Vietnam Institute for Advanced Study of Mathematics, 1 Dai Co Viet Street, Ha Noi, Viet Nam}
\email{lecongtrinh@qnu.edu.vn}
\author{Dư Thị-Hòa-Bình }
\address{Dư Thị-Hòa-Bình \\
Department of Mathematics, Quy Nhon University\\
170 An Duong Vuong, Quy Nhon, Binh Dinh}
\email{hoabinhcdsp@gmail.com}

\author{Nguyễn Trần-Đức }
\address{Nguyễn Trần-Đức \\
Department of Mathematics, Quy Nhon University\\
170 An Duong Vuong, Quy Nhon, Binh Dinh}
\email{nguyentranduc1995@gmail.com}

%    General info
%11E25  	Sums of squares and representations by other particular quadratic forms
%90C22  	Semidefinite programming
%13J30  	Real algebra [See also 12D15, 14Pxx]

\date{\today}

%\dedicatory{}


\begin{abstract} A number $\lambda \in \mathbb C $ is called an {\it eigenvalue} of the matrix polynomial $P(z)$ 
if there exists a nonzero vector $x \in \mathbb C^n$ such that $P(\lambda)x = 0$. 
Note that each finite eigenvalue of $P(z)$ is a zero of the   characteristic
polynomial $\det(P(z))$. In this paper we establish  some (upper and lower) bounds for eigenvalues of matrix polynomials based on the norm of their coefficient matrices and  compare these bounds to those given by N.J. Higham and F. Tisseur \cite{HT}, J. Maroulas and  P. Psarrakos \cite{MaPs}.
\end{abstract}

\maketitle
%%%%%%%%%%%%%%%%%%%%%%%%%%%%%%%%%%%%%%%%%%%%%%%%%%%%
%\tableofcontents
\section{Introduction}

Let $\mathbb C^{n\times n}$ be the set of all $n\times n$ matrices whose entries are in $\mathbb C$. For a \textit{matrix polynomial} we mean the matrix-valued function of a complex variable of the form 
\begin{equation}\label{mp}
P(z)= A_m z^m + \cdots + A_1 z + A_0, 
\end{equation}
where $A_i\in \mathbb C^{n\times n}$ for all $i=0,\cdots,m$. If $A_m\not =0$, $P(z)$ is called a matrix polynomial of \textit{degree} $m$. When $A_m=I$, the identity matrix in $C^{n\times n}$, the matrix polynomial $P(z)$ is called a \textit{monic}.

A number $\lambda \in \mathbb C$ is called an \textit{eigenvalue} of the matrix polynomial $P(z)$, if there exists a nonzero vector $x\in \mathbb C^n$ such that $P(\lambda)x=0$. Then the vector $x$ is called, as usual,  an \textit{eigenvector} of $P(z)$ associated  to the eigenvalue $\lambda$. Note that  each finite eigenvalue of $P(z)$ is a root  of  the \textit{characteristic polynomial} $\det(P(z))$.

The \textit{polynomial eigenvalue problem (PEP)} is to find an eigenvalue $\lambda$ and a non-zero vector $x\in \mathbb C^n$ such that $P(\lambda)x = 0$. For $m=1$, (PEP) is actually the \textit{generalized eigenvalue problem (GEP)} 
$$ Ax = \lambda Bx,$$
and, in addition, if $B=I$, we have the standard eigenvalue problem 
$$ Ax = \lambda x.$$ 
For $m=2$ we have the \textit{quadratic eigenvalue problem (QEP)}.

The theory of matrix polynomials was primarily devoted by two works, both of which are 
strongly motivated by the theory of \textit{vibrating systems}: one by Frazer, 
Duncan, and Collar in 1938 [FDC], and the other by P. Lancaster  in 1966 \cite{L}. 

(QEPs), and more generally (PEPs), play an important role in applications to science and engineering. We refer to \cite{TM} for a survey on applications of (QEP). Moreover, we refer to the book of  I. Gohberg, P. Lancaster and L. Rodman \cite{GLR} for a theory of matrix polynomials and their applications.

There are algorithms to solve (QEPs), see the works of Hamarling, Munro and Tisseur \cite[2013]{HMT}  and Zeng and Su \cite[2014]{ZS}. For (PEPs), there is some research on bounds of eigenvalues of matrix polynomials which were constructed in terms of the norms of coefficients of the given matrix polynomials. See, for example, the work of Higham  and Tisseur \cite[2003]{HT},    Maroulas and    Psarrakos \cite[1997]{MaPs}.

Computing eigenvalues of matrix polynomials (even computing eigenvalues of scalar matrices and finding roots of univariate  polynomials) is a hard problem.  There is an iterative method to compute these eigenvalues, see Simoncini and Perotti \cite[2006]{SP}. Moreover, when computing pseudospectra of matrix polynomials, which provide information about the global sensitivity of the eigenvalues, a particular region of the (possibly extended) complex plane must be identified that contains the eigenvalues of interest, and bounds clearly help to determine such region \cite{TH}. Therefore, it is useful  to find the location of these eigenvalues. 

  Note that, if $A_0$ is singular then $0$ is an eigenvalue of $P(z)$, and if $A_m$ is singular then $0$ is an eigenvalue of the matrix polynomial $z^mP(1/z)$. Therefore, to locate the eigenvalues of these matrix polynomials,  \textit{we always assume that $A_0$ and $A_m$ are non-singular}. 

The paper is organized as follows. In Section 2 we give bounds for matrix polynomials whose coefficients satisfy some special properties, in particular, we give a matrix version of Enestr\"om-Kakeya's theorem. In Section 3 we establish  matrix versions of some Cauchy's type theorems.  In particular,  we establish a matrix version of the theorem of Joyal, Labelle and Rahman (cf. \cite{JLR}, \cite[Theorem 2.14]{MR}) and some of its corollaries. Moreover, we give also a matrix version of  Datt and Govil's theorem     \cite[Theorem 1]{DG} and some other  bounds. Finally, we give some numerical experiments in Section 4.

\textbf{Notation}.  For a matrix $A\in \mathbb C^{n\times n}$, the notation   $A\geq 0$ means  \textit{"$A$ is positive semidefinite"}, i.e. for every vector $x\in \mathbb C^n$ we have $x^*Ax\geq 0$;   $A>0$ means  \textit{"$A$ is positive definite"}, i.e. $x^*Ax>0$ for every  nonzero vector $x\in \mathbb C^n$. 
For two matrices $A,B\in \mathbb C^{n\times n}$,  the notation  $A\geq B$ means  $A-B\geq 0$. 

\textit{Throughout this paper,  $\|\cdot \|$ denotes a subordinate matrix norm.}

\section{Enestr\"{o}m-Kakeya's theorem for matrix polynomials}
In this section we give upper and lower bounds for eigenvalues of some special matrix polynomials. First of all we consider matrix polynomials with a dominant property. 

\thr \label{dominance}  Let  $P(z)= A_0 + A_1 z+ \cdots+ A_m z^m$ be a matrix polynomial whose coefficients 
	 $A_i \in \mathbb{C}^{n \times n}$ satisfying the following  dominant property:
	 $$\| A_m\| > \|A_i\|, \forall~ i=0,\cdots,m-1. $$
Then all eigenvalues $\lambda$ of $P(z)$ locate in the open disk
 $$|\lambda|<1+\|A_m\| \|A_m^{-1}\|.$$
\ethr
In particular, for $n=1$, we obtain the following corollary of Cauchy's theorem (\cite[Theorem 27.2]{Ma}; see also   \cite[Theorem 2.2]{D}): \textit{Let  $p(z)= a_0 + a_1 z+ \cdots+ a_m z^m \in \mathbb C[z]$ such that $|a_m|>|a_i|$ for all $i=0,\cdots,m-1$. Then all the roots of $p(z)$ locate in the open disk $|z| <2$.}  The proof of this corollary uses the fact that when $n=1$ we have   $\|A_m\| \|A_m^{-1}\|=1$. 
\pf  Let $\lambda \in \mathbb C$ be an eigenvalue of $P(z)$ and $x\in \mathbb C^n$ a unit eigenvector of $P(z)$ associated to $\lambda$. \\
We have nothing to prove if $|\lambda| \leq 1$. Hence we may assume that $|\lambda| >1$. Then we have 
\begin{align*} 
\|P(\lambda)x\| & \geq |\lambda|^m\left[\|A_mx\| - \|\sum_{i=0}^{m-1}\dfrac{A_ix}{\lambda^{m-i}}\| \right]\\
&\geq |\lambda|^m\left[\|A_m^{-1}\|^{-1} - \sum_{i=0}^{m-1}\dfrac{\|A_i\|}{|\lambda|^{m-i}} \right]\\
&\geq |\lambda|^m\left[\|A_m^{-1}\|^{-1} - \sum_{i=0}^{m-1}\dfrac{\|A_m\|}{|\lambda|^{m-i}} \right]\\
&= |\lambda|^m\|A_m^{-1}\|^{-1}\left[1 - \|A_m\|\|A_m^{-1}\|\sum_{i=1}^{m}\dfrac{1}{|\lambda|^{i}} \right]\\
&> |\lambda|^m\|A_m^{-1}\|^{-1}\left[1 - \|A_m\|\|A_m^{-1}\|\sum_{i=1}^{\infty}\dfrac{1}{|\lambda|^{i}} \right]\\
&=|\lambda|^m\|A_m^{-1}\|^{-1}\left[1 - \dfrac{\|A_m\|\|A_m^{-1}\|}{|\lambda| -1} \right]\\
&= \dfrac{|\lambda|^m\|A_m^{-1}\|^{-1}}{|\lambda| -1} \left(|\lambda| -1 -\|A_m\|\|A_m^{-1}\|\right).
\end{align*}
Hence, if $|\lambda| \geq 1+\|A_m\|\|A_m^{-1}\|$ we have  $\|P(\lambda)x\|>0$, a contradiction. It follows that $|\lambda| <  1+\|A_m\|\|A_m^{-1}\|$, which completes the proof.
\epf 

The  following theorem of Enestr\"{o}m and Kakeya  is well-known.

\thr[{\cite[Corollary 3]{SiSh}}] \label{dang1}
Let $p(z)$ be a polynomial in one variable given by 
$$p(z)= a_0 + a_1 z+ \cdots+ a_mz^m, ~a_i\in \mathbb R, \forall i=1,\cdots,m.$$ 
Suppose that 
$$a_m \geq a_{m-1} \geq \cdots \geq a_0 \geq 0;~ a_m>0.$$
If   $z\in \mathbb C$ is a root of $p(z)$ then  $\dfrac{a_0}{2a_m}\leq |z|\leq 1$.
 \ethr
 
 A matrix version of Theorem \ref{dang1} is given as follows.
 
 \thr\label{thrEK1} 
	Let  $P(z)= A_0 + A_1 z+ \cdots+ A_m z^m$ be a matrix polynomial whose coefficients 
  $A_i\in \mathbb C^{n\times n}$  satisfying 
	$$A_m \geq A_{m-1} \geq \cdots \geq A_0 \geq 0; ~A_m>0.$$
	Then  each eigenvalue  $\lambda$ of $P(z)$ satisfies   
	$$ \dfrac{\lambda_{\min}(A_0)}{2\lambda_{\max}(A_m)} \leq \left| \lambda \right|  \leq 1,$$
	where $\lambda_{\min}(A_0)$ denotes the smallest eigenvalue of $A_0$ and   
	 $\lambda_{\max}(A_m)$ the largest eigenvalue of $A_m$.
		\ethr
	\pf A proof for the upper bound of $|\lambda|$ in this theorem was given by G. Dirr and H. K. Wimmer \cite[Theorem 2.1]{DW} (see also in \cite[Theorem 5.1]{SW}). Now we give a proof for the lower bound. 	
	
	Firstly we observe that for a matrix $A\in \mathbb C^{n\times n}$, its smallest  eigenvalue $\lambda_{min}(A)$ and its largest eigenvalue   $\lambda_{max}(A)$  belong to the set 
	$$\{x^*Ax | x\in \mathbb C^n, \|x\|=1\},$$
	which is the \textit{standard numerical range of $A$}. 
	Hence  for a unit vector $x\in \mathbb C^n$, we always have 
\begin{equation}\label{pt:minmax}
	   \lambda_{min}(A) \leq x^*Ax \leq \lambda_{max}(A).
\end{equation}
Let $\lambda \in \mathbb C$ be an eigenvalue of $P(z)$, and $u\in \mathbb C^n, \|u\|=1$ an eigenvector of $P(z)$ associated  to $\lambda$. Consider the polynomial 
$$P_u(z): = u^*P(z)u = \sum_{i=0}^{m}(u^*A_iu)z^i.$$
	Note that   $\lambda $ is a root of  $P_u(z)$. Moreover, the hypothesis on the relation of $A_i$'s implies that 
	$$u^*A_mu \geq u^*A_{m-1}u \geq \cdots \geq u^*A_0u \geq 0, u^*A_mu >0, $$
that is, the polynomial $P_u(z)$ satisfies the conditions given in Theorem \ref{dang1}. 	
	Applying this theorem for $P_u(z)$ we obtain   	
	$$  \dfrac{u^*A_0u}{2u^*A_mu}\leq |\lambda|. $$
	Then the required lower bound for $|\lambda|$ follows from (\ref{pt:minmax}).
\epf

 By applying Theorem \ref{thrEK1} for the matrix polynomial $z^nP(\frac{1}{z})$  we obtain the following dual version  of Theorem \ref{thrEK1}.
 \thr\label{thrEK1-dual} 
	Let  $P(z)= A_0 + A_1 z+ \cdots+ A_m z^m$ be a matrix polynomial whose coefficients 
  $A_i\in \mathbb C^{n\times n}$  satisfying 
	$$A_0 \geq A_1 \geq \cdots \geq A_m > 0.$$
	Then  each eigenvalue  $\lambda$ of $P(z)$ satisfied     
	$ \left| \lambda \right|  \geq 1.$
	\ethr

  We have also the following version of Enestr\"{o}m-Kakeya's theorem for polynomials. 
\thr[Enestr\"{o}m-Kakeya's theorem, Version 2, \cite{BE}] \label{dang2}
Let  $p(z)= a_0 + a_1 z+ \cdots+ a_mz^m$ be a polynomial whose coefficients  $a_i, i=0,\cdots, m$ are positive real numbers. Denote
$$\alpha:= \min_{0 \leq i \leq m-1} \left\{ \frac{a_i}{a_{i+1}}\right\}, \beta:= \max_{0 \leq i \leq m-1} \left\{ \frac{a_i}{a_{i+1}}\right\}.$$
Then each root $z\in \mathbb C $ of $p(z)$ satisfies the following inequalities 
$$\alpha \leq \left| z \right| \leq \beta.$$
\ethr
Using the same method as given in the proof of Theorem   \ref{thrEK1}, applying Theorem \ref{dang2}, we obtain  the following bounds for eigenvalues of matrix polynomials  whose coefficients are positive definite.   
	
\thr	\label{thrEK2}
	Let  $P(z)= A_0 + A_1 z+ \cdots+ A_m z^m$ be a matrix polynomial whose coefficients 
	 $A_i \in \mathbb{C}^{n \times n}$ are positive definite.    If $\lambda\in \mathbb C$ is an eigenvalue of  $P(z)$, then
	$$\min_{i=0,\cdots,m-1}\left\{\frac{\lambda_{min}(A_i)}{\lambda_{max}(A_{i+1})} \right\} \leq \left| \lambda\right| \leq \max_{i=0,\cdots,m-1}\left\{\frac{\lambda_{max}(A_i)}{\lambda_{min}(A_{i+1})} \right\}.$$
	\ethr

\section{Cauchy type theorems for  matrix polynomials }
In this section we establish some Cauchy type theorems for matrix polynomials of the form  $P(z)= A_0 + A_1 z+ \cdots+ A_m z^m$ with $A_m$ and $A_0$ non-singular. We should observe that the set of eigenvalues of $P(z)$ coincides to that of the monic matrix polynomial    
$$A_m^{-1}P(z)= (A_m^{-1}A_0) + (A_m^{-1}A_1) z+ \cdots+ I z^m.$$
Therefore, because of the complexity in practice, \textit{we concentrate to consider in this section the bounds for monic matrix polynomials}. 

Firstly we state the Cauchy's theorem for monic matrix polynomials.
  
\thr[{Cauchy, \cite[Lemma 3.1]{HT}}]\label{Cauchy-Matrix}
Let  $P(z)= A_0 + A_1 z+ \cdots+ I z^m$ be a monic matrix polynomial. Let  $r$ resp. $R$ be   the positive root of the  polynomial
$$ h(z)= z^m + z^{m-1} \left \| A_{m-1} \right \| + \cdots+ z\left \| A_1 \right \| - \left \| A_0^{-1} \right \|^{-1}
$$ 
resp. 
$$ g(z)= z^m - z^{m-1} \left \| A_{m-1} \right \| - \cdots - \left \| A_0 \right \|. $$
	Then each eigenvalue  $\lambda$ of $P(z)$ satisfies
	$$r \leq \left | \lambda \right | \leq R.$$
\ethr
Now we give some Cauchy type theorem for monic matrix polynomials.

\thr  \label{Cauchytype1} Let  $P(z)= A_0 + A_1 z+ \cdots+ I z^m$. Denote $$M:= \displaystyle \max_{i=0,\cdots,m-1}  \|A_i\|.$$
Then all eigenvalues of $P(z)$ are contained  in the closed disk 
$$K(0,r_1) :=\{ z\in \mathbb C |~~ \left|z\right| \leq r_1\}, $$ where $r_1:=\max\{1,\delta\}$ and $\delta \not =1$ is the positive root of the equation 
$$ z^{m+1} - (1+M)z^m + M =0. $$
\ethr
In particular, for $n=1$ we obtain a Cauchy type theorem for polynomials \cite[Theorem 3.2]{D}. 

\pf Let $\lambda \in \mathbb C$ be an eigenvalue of $P(z)$ and $x\in \mathbb C^n$ a unit eigenvector of $P(z)$ associated to $\lambda$. \\
The conclusion is clear if $|\lambda| \leq 1$. Therefore we may assume that $|\lambda| >1$. Then we have 
\begin{align} 
\|P(\lambda)x\| & \geq \left[\|Ix\||\lambda|^m - \|\sum_{i=0}^{m-1} A_ix \lambda^{i}\| \right] \nonumber\\
&\geq  \left[|\lambda|^m  - \sum_{i=0}^{m-1}\|A_i\| \|A_m^{-1}\| \lambda^{i} \right]\label{line2}\\
&\geq   \left[|\lambda|^m  - M\sum_{i=0}^{m-1} \lambda^i \right]\label{line3} \\
& = \left[|\lambda|^m  - M \dfrac{|\lambda|^m - 1}{|\lambda| -1} \right]\nonumber\\
&=\dfrac{1}{|\lambda| -1}\left(|\lambda|^{m+1}-(1+M)|\lambda|^{m}+M\right).\nonumber
\end{align}
In the lines above, from (\ref{line2}) to (\ref{line3}) we use the definition of $M$.\\
Note that the polynomial $f(z):=z^{m+1}-(1+M)z^m + M$ has exactly two positive real roots $1$ and $\delta\not =1$  by the Descartes' rule of signs, and $f(0)>0$. It follows that 
$$ |f(z)| > 0 \mbox{ for all } z >  \max\{\delta, 1\}. $$
Hence for $|\lambda| >  r_1$ we have $\|P(\lambda)x\| >0$, a contradiction. This completes the proof.
\epf

\coro  \label{Cauchytype2} Let  $P(z)= A_0 + A_1 z+ \cdots+ I z^m$ be a monic matrix polynomial. Denote $$\widetilde{M}:=\displaystyle \max_{i=0,\cdots,m}  \|A_{m-i}-A_{m-i-1}\| \quad (A_m=I \mbox{ and } A_{-1}=0).$$
Then all eigenvalues of $P(z)$ are contained  in the closed disk $K(0,r_2)$, where $r_2:=\max\{1,\delta\}$ and  $\delta \not =1$ is the positive root of the equation 
$$ z^{m+2} - (1+\widetilde{M})z^{m+1} + \widetilde{M} =0. $$
\ecoro
In particular, for $n=1$ we obtain   \cite[Theorem 3.3]{D}. 
\pf Consider the matrix polynomial 
$$Q(z):=(1-z)P(z)=-Iz^{m+1}+\sum_{i=0}^m(A_{m-i}-A_{m-i-1})z^{m-i}. $$
Applying Theorem \ref{Cauchytype1} for the polynomial $Q(z)$, observing that each eigenvalue of $P(z)$ is also an eigenvalue of $Q(z)$, we obtain the required result.
\epf
	
\thr  \label{Cauchytype3} Let  $P(z)= A_0 + A_1 z+ \cdots+ I z^m$. Then all eigenvalues of $P(z)$ are contained  in the open disk 
$$K^o(0,r_3):=\{ z\in \mathbb C |~~ \left|z\right| < r_3\},$$ where  $r_3:=1+{M}$ and ${M}$ is defined as in Theorem \ref{Cauchytype1}.
\ethr	
In particular, for $n=1$ we obtain   another Cauchy's theorem for polynomials \cite[Theorem (27,2)]{Ma}.
\pf Let $\lambda \in \mathbb C$ be an eigenvalue of $P(z)$ and $x\in \mathbb C^n$ a unit eigenvector of $P(z)$ associated to $\lambda$. \\
As above, we may assume that $|\lambda| >1$.
Then we have 
\begin{align*} 
\|P(\lambda)x\| & \geq |\lambda|^m\left[\|Ix\| - \|\sum_{i=0}^{m-1}\dfrac{A_ix}{\lambda^{m-i}}\| \right]\\
&\geq |\lambda|^m\left[1 - \sum_{i=0}^{m-1}\dfrac{\|A_i\|}{|\lambda|^{m-i}} \right]\\
&\geq |\lambda|^m\left[1 - M\sum_{i=1}^{m}\dfrac{1}{|\lambda|^{i}} \right]\\
&> |\lambda|^m\left[1 - M\sum_{i=1}^{\infty}\dfrac{1}{|\lambda|^{i}}\right]\\
&=|\lambda|^m\left[1 - \dfrac{M}{|\lambda| -1} \right]\\
&= \dfrac{|\lambda|^m}{|\lambda| -1} \left(|\lambda| -1 -M\right).
\end{align*}
Then, for $|\lambda| \geq 1+M$ we have $\|P(\lambda)x\| > 0$, a contradiction. Thus $|\lambda| < 1+M$.
\epf

\coro  \label{Cauchytype4} Let  $P(z)= A_0 + A_1 z+ \cdots+ I z^m$. Then all eigenvalues of $P(z)$ are contained  in the open disk 
$K^o(0,r_4)$,   where $r_4:=1+\widetilde{M}$ and  $\widetilde{M}$ is defined as in Corollary \ref{Cauchytype2}. 
\ecoro	
In particular, for $n=1$ we obtain   \cite[Theorem 3.4]{D}.
\pf Consider the matrix polynomial 
$$Q(z):=(1-z)P(z)=-Iz^{m+1}+\sum_{i=0}^m(A_{m-i}-A_{m-i-1})z^{m-i}. $$
Since each eigenvalue of $P(z)$ is also an eigenvalue of   $Q(z)$, applying Theorem \ref{Cauchytype3} for $Q(z)$ we have the conclusion.  
\epf 

Next we give a matrix version of the theorem of Joyal, Labelle and Rahman, cf. \cite{JLR}, \cite[Theorem 2.14]{MR}.
\thr \label{Cauchytype5} Let  $P(z)= A_0 + A_1 z+ \ldots+ A_{m-1}z^{m-1}+ I  z^m$ be a monic matrix polynomial. Denote
$$\alpha:= \displaystyle \max_{i={0,\cdots,m-2}} \left\| A_i \right\|.$$ 
Then each eigenvalue  $\lambda$ of  $P(z)$ is estimated by
$$ \left| \lambda \right| \leq \frac{1}{2} \left\lbrace 1+ \left\| A_{m-1} \right\| + \left[\left( 1- \left\| A_{m-1} \right\|\right)^2 +4 \alpha  \right] ^\frac{1}{2} \right\rbrace. $$
\ethr

\pf Let $\lambda \in \mathbb C$ be an eigenvalue of $P(z)$ and $x\in \mathbb C^n$ a unit eigenvector of $P(z)$ associated to $\lambda$. \\
By contradiction, assume 
$$ \left|\lambda \right| > \frac{1}{2} \left\lbrace 1+ \left\| A_{m-1} \right\| + \left[\left( 1- \left\| A_{m-1} \right\|\right)^2 +4 \alpha \right] ^\frac{1}{2} \right\rbrace.$$
It follows that 
\begin{equation}\label{equ-Cauchy5}
\left( \left| \lambda \right| -1 \right) \left( \left| \lambda \right| - \left\| A_{m-1}\right\|  \right) - \alpha >0.
\end{equation}
Multiplying  (\ref{equ-Cauchy5}) by $|\lambda|^{m-1}$ and then dividing by $|\lambda| -1$, we obtain
$$ \left| \lambda \right| ^m - \left\| A_{m-1}\right\| \lambda^{m-1} - \alpha \frac{\left| \lambda \right| ^{m-1}}{\left|  \lambda \right| -1} >0. $$
However,
\begin{align*}
\alpha\frac{\left| \lambda \right| ^{m-1}}{\left|  \lambda \right| -1} & >  \alpha \frac{\left| \lambda \right| ^{m-1} -1}{\left|  \lambda \right| -1}  =\alpha(1 + \left| \lambda \right| +\cdots +\left|  \lambda \right|^{m-2}) \\
& \geq \|(A_0+A_1\lambda  +\cdots + A_{m-2}\lambda^{m-2})x\|.
\end{align*} 
On the other hand, 
$$  \left| \lambda \right| ^m - \left\| A_{m-1}\right\| \lambda^{m-1} \leq \|(I\cdot \lambda^m + A_{m-1}\lambda^{m-1})x\|.$$
It follows that
\begin{align*}
 0&< \left| \lambda \right| ^m - \left\| A_{m-1}\right\| \lambda^{m-1} - \alpha \frac{\left| \lambda \right| ^{m-1}}{\left|  \lambda \right| -1} \\
 & < \|(I\cdot \lambda^m + A_{m-1}\lambda^{m-1})x\| - \|(A_0+A_1\lambda  +\cdots + A_{m-2}\lambda^{m-2})x\| \\
 & \leq \|(A_0+A_1\lambda  +\cdots + A_{m-2}\lambda^{m-2})x+(A_{m-1}\lambda^{m-1}+I\cdot \lambda^m)x\|=\|P(\lambda)x\|,   
 \end{align*}
 a contradiction. Thus
$$ \lambda \leq \frac{1}{2} \left\lbrace 1+ \left\| A_{m-1} \right\| + \left[\left( 1- \left\| A_{m-1} \right\|\right)^2 +4 \alpha \right] ^\frac{1}{2} \right\rbrace. $$
\epf 

By applying Theorem  \ref{Cauchytype5} for the monic matrix polynomial $z^mP(\frac{1}{z})$ we obtain  the following lower bound for eigenvalues of $P(z)$.

\coro \label{coro1}  Let  $P(z)= A_0 + A_1 z+ \cdots+ A_{m-1}z^{m-1}+I z^m$.  Denote $L_i:=A_0^{-1}A_i$ ($i=1,\ldots,m-1)$,   $L_m=A_0^{-1}$, and
$$\beta:= \displaystyle \max_{i={2,\cdots,m}} \left\| L_i \right\|. $$ 
Then for each eigenvalue  $\lambda$ of  $P(z)$ we have
$$ \left| \lambda \right| \geq \dfrac{2}{ 1+ \left\| L_1 \right\| + \left[\left( 1- \left\| L_1\right\|\right)^2 +4 \beta  \right] ^\frac{1}{2}}. $$
\ecoro 
By applying Theorem  \ref{Cauchytype5} for the matrix polynomial $(1-z)P(z)$ we obtain 
\coro \label{coro2}  Let  $P(z)= A_0 + A_1 z+ \cdots+ A_{m-1}z^{m-1}+I z^m$. Denote
$$\gamma:= \displaystyle \max_{i={1,\cdots,m}} \left\| A_{m-i}-A_{m-i-1}\right\| \quad \quad \quad  (A_{-1}=0).$$ 
Then   each eigenvalue  $\lambda$ of  $P(z)$ satisfies
$$ \left| \lambda \right| \leq \frac{1}{2} \left\lbrace 1+ \left\| I - A_{m-1}\right\| + \left[\left( 1- \left\| I - A_{m-1} \right\|\right)^2 +4 \gamma  \right] ^\frac{1}{2} \right\rbrace. $$
\ecoro 
Similarly, Corollary \ref{coro2} yields  the following lower bound. 

\coro \label{coro2'}  Let  $P(z)= A_0 + A_1 z+ \ldots+ A_{m-1}z^{m-1}+ I z^m$.  Denote
$$\gamma':= \displaystyle \max_{i={1,\cdots,m}} \left\| L_i-L_{i+1}\right\| \quad \quad \quad (L_{m+1}=0).$$ 
Then  each eigenvalue  $\lambda$ of  $P(z)$ satisfies 
$$ \left| \lambda \right| \geq \dfrac{2}{1+ \left\| I-L_1\right\| + \left[\left( 1- \left\| I-L_1 \right\|\right)^2 +4 \gamma'  \right] ^\frac{1}{2}}. $$
\ecoro

By applying Theorem \ref{Cauchytype5} for the matrix polynomial $(I z - A_{m-1})P(z)$ we obtain 
\coro \label{coro3}  Let  $P(z)= A_0 + A_1 z+ \ldots+ A_{m-1}z^{m-1}+ I  z^m$. Denote
$$\delta:= \displaystyle \max_{i={0,\cdots,m-1}} \left\| A_{m-1}A_i-A_{i-1} \right\| \quad \quad \quad (A_{-1}=0).$$ 
Then   each eigenvalue  $\lambda$ of  $P(z)$ satisfies
$$ \left| \lambda \right| \leq \frac{1}{2} (1+\sqrt{1+4\delta}). $$
\ecoro

Corollary \ref{coro3} yields  the following lower bound.
\coro \label{coro3"}
Let Let  $P(z)= A_0 + A_1 z+ \ldots+ A_{m-1}z^{m-1}+ I  z^m$. Denote
$$\delta':= \displaystyle \max_{i={1,\cdots,m}} \left\| L_1L_i-L_{i+1}\right\| \quad \quad \quad (L_{m+1}=0).$$ 
Then each eigenvalue  $\lambda$ of  $P(z)$ is bounded below by
$$ \left| \lambda \right| \geq \dfrac{2}{1+\sqrt{1+4\delta'}}. $$
\ecoro

By applying Theorem \ref{Cauchytype5} for the matrix polynomial $(I\cdot z + I - A_{m-1})P(z)$ we obtain 
\coro \label{coro4}  Let  $P(z)= A_0 + A_1 z+ \ldots+ A_{m-1}z^{m-1}+ I z^m$. Denote
$$\epsilon:= \displaystyle \max_{i={0,\cdots,m-1}} \left\| (I-A_{m-1})A_i+A_{i-1} \right\| \quad \quad \quad (A_{-1}=0).$$ 
Then  each eigenvalue  $\lambda$ of  $P(z)$ is estimated by 
$$ \left| \lambda \right| \leq 1+\sqrt{\epsilon}. $$
\ecoro 

The following lower bound is  obtained by applying Corollary \ref{coro4} for the matrix polynomial $z^mP(\frac{1}{z})$.

\coro \label{coro4"}
 Let  $P(z)= A_0 + A_1 z+ \ldots+ A_{m-1}z^{m-1}+ I z^m$. Denote
$$\epsilon':= \displaystyle \max_{i={1,\cdots,m}} \left\| (I-L_1)L_i+L_{i+1}\right\| \quad \quad \quad (L_{m+1}=0).$$ 
Then  each eigenvalue  $\lambda$ of  $P(z)$ bounded below by 
$$ \left| \lambda \right| \geq \dfrac{1}{1+\sqrt{\epsilon'}}. $$
\ecoro

Next we give the matrix version of the theorem of Datt and Govil     \cite[Theorem 1]{DG}.

\thr \label{Cauchytype6}Let  $P(z)= A_0 + A_1 z+ \cdots+ I  z^m$ be a monic matrix polynomial. Denote 
$$M = \displaystyle\max_{i=0,\cdots, m-1} \left\|A_i \right\|. $$
 Then  each eigenvalue $\lambda$ of $P(z)$ satisfies 
$$\frac{\left\|A_0^{-1} \right\|^{-1}}{2(1+M)^{m-1}(Mm+1)} \leq \left|\lambda \right| \leq 1 + \lambda_0M,$$
where  $\lambda_0$ is a root of the equation $x = 1 - \frac{1}{(Mx+1)^m}$ in the interval $(0,1)$.
 \ethr
  \pf 
 Let $\lambda \in \mathbb C$ be an eigenvalue of $P(z)$ and $x\in \mathbb C^n$ a unit eigenvector of $P(z)$ associated to $\lambda$. 
 
 First we prove the upper bound for $\left|\lambda \right|$. We consider two cases:\\ 
\textbf{The first case:}  $mM \leq 1$. In this case,  if   $\left| \lambda\right| > 1$, we have
$$\left\| P(\lambda)x \right\|  \geq \left| \lambda\right| ^{m} - mM\left| \lambda\right| ^{m-1} \geq \left| \lambda\right| ^{m} - \left|\lambda \right|^{m-1} > 0, \text{ a contradiction.}$$
It follows that  $\left| \lambda\right| \leq 1 \leq  1 + \lambda_0M $ for all $\lambda_0 \in (0,1).$\\
\textbf{The second case:}  $mM>1$. In this case the equation  $x = 1 - \frac{1}{(Mx+1)^m}$ has  a unique root    $\lambda_0\in (0,1)$ \cite[Lemma 2]{DG}. Moreover, we have 
$$\left\|P(\lambda)x \right\| \geq \left|\lambda \right|^{m}  - M\sum_{j=0}^{m-1}\left|\lambda \right|^{j} = \left|\lambda \right|^{m} - M\frac{\left|\lambda \right|^{m} - 1  }{\left|\lambda \right|-1} .$$
If $\left|\lambda \right| > 1 + M\lambda_0$, we can write  $\left|\lambda \right| = 1 + M\alpha$ with $\alpha > \lambda_0$. Then $\alpha > 1- \frac{1}{(M\alpha +1)^m}$. It follows that 
$$\left\| P(\lambda)x \right\|  \geq (1+M\alpha)^m - \frac{(1+M\alpha)^m-1}{\alpha} >0,$$
a contradiction.   Thus $\left|\lambda \right| \leq 1 + M\lambda_0.$\\
Now we prove the lower bound for $|\lambda|$.  By contradiction, assume  $\left|\lambda \right| < \frac{\left\|A_0^{-1} \right\|^{-1}}{2(1+M)^{m-1}(Mm+1)}.$
Let us consider the matrix polynomial  $G(z) := (1-z)P(z).$\\
We have $$G(z) = A_0 + \sum_{i=1}^{m}(A_i - A_{i-1})z^i + I  z^m - A_{m-1}z^m - I  z^{m+1}=: A_0 + H(z).$$
Denote $R:=1+M$. Then for $|z|=R$, we have 
$$\begin{aligned}
\max_{\left|z \right| = R } \left\|H(z)x \right\| &\leq R^{m+1} + R^{m} + \left\| A_{m-1}\right\| R^m + \sum_{i=1}^{m-1}\left\|A_i - A_{i-1} \right\|R^i\\
&\leq R^m\left[R + 1 + M + 2(m-1)M \right]\\
&= 2(1+M)^m(mM+1).\\   
\end{aligned} $$
It follows from the \textit{maximal module principle} that for $|z|\leq R$ we have 
$$ \left\|H(z)x \right\| \leq 2(1+M)^m(mM+1).$$
Then for  $\left|\lambda \right| < \frac{\left\|A_0^{-1} \right\|^{-1}}{2(1+M)^{m-1}(Mm+1)} < R$ we have
$$\begin{aligned}
\left\| G(\lambda)x \right\| = \left\| A_0x + H(\lambda)x \right\| &\geq \left\| A_0^{-1}\right\| ^{-1} - \left\| H(\lambda)x\right\| \\
&\geq \left\| A_0^{-1}\right\| ^{-1} - \frac{\left| \lambda\right| }{1 + M}\max_{\left|\lambda \right| \leq 1 + M } \left\|H(\lambda)x \right\|\\
&\geq \left\| A_0^{-1}\right\| ^{-1} - 2(1+M)^{m-1}(mM+1)\left| \lambda\right| > 0 , 
\end{aligned} $$
a contradiction. Therefore  $$\frac{\left\|A_0^{-1} \right\|^{-1}}{2(1+M)^{m-1}(Mm+1)} \leq \left|\lambda \right|.$$
  \epf 
  
 If we do not wish to look for a root in the interval $(0,1)$ of the equation $x = 1 - \frac{1}{(Mx+1)^m}$, we use the following upper bound.
 \coro \label{coro5'} Let  $P(z)= A_0 + A_1 z+ \cdots+ I  z^m$ be a monic matrix polynomial.
 Then   each eigenvalue $\lambda$ of $P(z)$ satisfies
$$\frac{\left\|A_0^{-1} \right\|^{-1}}{2(1+M)^{m-1}(Mm+1)} \leq \left|\lambda \right| < 1 + \Big(1-\dfrac{1}{(1+M)^m}\Big)M.$$
  \ecoro
  \pf The proof follows from Theorem  \ref{Cauchytype6} and  the fact that for a root $\lambda_0$ of the equation $x = 1 - \frac{1}{(Mx+1)^m}$ in the interval $(0,1)$, we have always $ \lambda_0 < 1-\dfrac{1}{(1+M)^m}.$  
  \epf 
	
Next we give some other bounds for the magnitude of eigenvalues of monic matrix polynomials.

\thr \label{other1} Let  $P(z)= A_0 + A_1 z+ \cdots+ I   z^m$ be a monic matrix polynomial.  Denote 
$$ M:=\displaystyle \max_{i=0,\cdots,m-1} \left \| A_i \right \|, ~~M':= \displaystyle \max_{i=1,\cdots,m} \left \| A_i \right \|.$$
Then each eigenvalue $\lambda $  of $P(z)$ satisfies 
$$ \frac{\left \| A_0 ^{-1} \right \| ^{-1}}{\left \| A_0 ^{-1} \right \| ^{-1} + M' } < \left | \lambda \right | < 1 + M. $$
\rm In particular, for $n=1$ we obtain \cite[Theorem 2.2]{MR}. 
\ethr
\pf The upper bound is the one  obtained in Theorem   \ref{Cauchytype3}. Now we prove the lower bound. \\
Let $\lambda \in \mathbb C$ be an eigenvalue of $P(z)$ and $x\in \mathbb C^n$ a unit eigenvector of $P(z)$ associated to $\lambda$. 
If  $\left | \lambda \right | \leq \frac{\left \| A_0 ^{-1} \right \| ^{-1}}{\left \| A_0 ^{-1} \right \| ^{-1} + M'}$, we have 
 \begin{align*} 
\left \| P(\lambda)x \right \| &\geq \left \| A_0 ^{-1} \right \| ^{-1} - \sum_{i=1}^{m} \left | \lambda \right | ^i \left \| A_i \right \| \\
& \geq \left \| A_0 ^{-1} \right \| ^{-1} - M' \sum_{i=1}^{m} \left | \lambda \right | ^i  \\
& > \left \| A_0 ^{-1} \right \| ^{-1} - M' \frac{\left | \lambda \right |}{1-\left | \lambda \right |} \\
& = \frac{\left \| A_0 ^{-1} \right \| ^{-1} \left( 1-\left | \lambda \right |\right)  - M' \left | \lambda \right |}{1-\left | \lambda \right |} \geq  0, \text{ ~ a contradiction}.
\end{align*}
It follows that  $\left | \lambda \right | > \frac{\left \| A_0 ^{-1} \right \| ^{-1}}{\left \| A_0 ^{-1} \right \| ^{-1} + M'}.$   This completes the proof.
\epf 
More generally, we have the following bounds.
\thr \label{other2} Let  $P(z)= A_0 + A_1 z+ \cdots+ I  z^m$ be a monic matrix polynomial. 
  Let  $p,q >1 $ such that $ \frac{1}{p} + \frac{1}{q} = 1$. Denote 
  $$M_p := \left( \displaystyle\sum_{i=0}^{m-1} \left\| A_i \right\|^p \right)^\frac{1}{p}, ~ M'_p := \left( \displaystyle\sum_{i=1}^{m} \left\| A_i \right\|^p \right)^\frac{1}{p}.$$
Then each eigenvalue  $\lambda$  of $P(z)$ satisfies
$$ \left [\frac{\left \| A_0 ^{-1} \right \| ^{-q}}{(M'_p)^q + \left \| A_0 ^{-1} \right \| ^{-q}} \right ]^\frac{1}{q} < \left | \lambda \right | < (1 + M_p^q)^\frac{1}{q}. $$
\ethr
 In particular, for $n=1$ we obtain \cite[Theorem 2.4]{MR}. Moreover, letting $p$ tend to infinity (then $q$ tends to $1$), we obtain Theorem \ref{other1}.
\pf 
Let $\lambda \in \mathbb C$ be an eigenvalue of $P(z)$ and $x\in \mathbb C^n$ a unit eigenvector of $P(z)$ associated to $\lambda$. \\
If $ \left | \lambda \right | \geq (1+M_p^q)^\frac{1}{q} $, we have
\begin{align} 
	\left \| P(\lambda)x \right \| &\geq   \left| \lambda \right|^m - \sum_{i=0}^{m-1} \left\| A_i\right\| \left| \lambda \right|^i\label{line4.1}\\
	& \geq    \left| \lambda \right|^m - \left(\sum_{i=0}^{m-1} \left\| A_i\right\|^p \right)^\frac{1}{p} \left(\sum_{i=0}^{m-1} \left| \lambda \right|^{iq} \right)^\frac{1}{q}\label{line4.2} \\
	& =    \left| \lambda \right|^m \left[1- \frac{M_p}{\left| \lambda \right|^m} \left(\sum_{i=0}^{m-1} \left| \lambda \right|^{iq} \right)^\frac{1}{q}   \right] \nonumber\\
	& =  \left| \lambda \right|^m \left[1-  M_p \left(\sum_{i=0}^{m-1} \left| \lambda \right|^{(i-m)q} \right)^\frac{1}{q}   \right]\nonumber \\
	& >   \left| \lambda \right|^m\left[1-  M_p  \left(\sum_{i=1}^{\infty} \left| \lambda \right|^{-iq} \right)^\frac{1}{q}   \right] \nonumber \\
	& = \left| \lambda \right|^m \left[1-  M_p  \frac{1}{\left( \left|\lambda\right| ^q - 1 \right)^\frac{1}{q} }  \right] \geq 0, \text{ a contradiction}.\nonumber
\end{align}
In the lines above, from (\ref{line4.1}) to (\ref{line4.2}) we use the well-known H\"{o}lder's inequality. \\
It follows that  $ \left | \lambda \right | < (1 + M_p^q)^\frac{1}{q}. $\\
Similarly  we have  $\left | \lambda \right | > \left [\frac{\left \| A_0 ^{-1} \right \| ^{-q}}{(M'_p)^q + \left \| A_0 ^{-1} \right \| ^{-q}} \right ]^\frac{1}{q}. $ This completes the proof.
\epf

\section{Numerical experiments}
We have  already established several  estimations for eigenvalues of matrix polynomials. It is in general not possible to compare the sharpness of these bounds. We can only compare them in some special cases by  numerical examples. In order to get a good comparison throughout practical examples, we use random data in each example.  Moreover, we compare the sharpness of our bounds and those given by N.J. Higham and F. Tisseur \cite{HT}, J. Maroulas and  P. Psarrakos \cite{MaPs}. We compute  and compare the bounds for two cases of the matrix coefficients: One with arbitrary random matrix coefficients, and the other one with symmetric matrix coefficients.  The experiments were performed using the open source software OCTAVE (version 4.4.0).

\begin{example} \rm 

Consider a $5\times 5$ monic matrix polynomial $P(z)$ of degree $m=9$ whose coefficient matrices are given by
$$A_i=10^{i-3}rand(5), ~ i=0,\ldots, 8, $$
where $rand(5)$ denotes a $5 \times 5$ random matrix from the normal $(0,1)$ distribution.

The upper bounds obtained by Higham and Tisseur \cite{HT}  are given in Table \ref{higham-upperbound}, while our new upper bounds are given in Table \ref{upperbound-new}.  
\begin{center}
\begin{table}[h!]
\begin{tabular}{|l|l|l|} 
\hline 
\textbf{Lemmas \quad } & \textbf{Values} & \textbf{Comments \quad \quad \quad \quad \quad \quad \quad \quad \quad\quad   }\\
\hline 
2.3 ~(2.2) & 3.5422 $\times 10^5$ & $\infty$-norm based\\
\hline
2.3 ~ (2.3) & 2.4987 $\times 10^5$ &2-norm based \\
\hline
2.5~ (2.13)& 3.3493 $\times 10^5 $ & $\infty$-norm based\\
\hline
2.6 ~ (2.14)&  3.4651$\times 10^5 $& Ostrowski, $\beta=3/4$\\
\hline
2.11~ (2.18) &  2.4907 $\times 10^5$ & 2-norm based\\
\hline
3.1 & 2.4827 $\times 10^5$ & Cauchy's theorem applied for $P$, 2-norm \\
\hline
3.1 & 2.4827 $\times 10^5$ & Cauchy's theorem applied for $P_U$, 2-norm \\
\hline
4.1& 4.9654 $\times 10^5$ & 2-norm based\\
\hline
\end{tabular}
\vspace{0.1cm}
\caption{Higham and Tisseur's upper bounds}  \label{higham-upperbound}
\end{table}

\begin{table}[h!]
\begin{tabular}{|l|l|l|}
\hline 
\textbf{Theorems/Corollaries} & \textbf{Values} & \textbf{Comments \quad \quad \quad \quad \quad \quad \quad  \quad }\\
\hline 
3.2, 3.2.1, 3.3, 3.3.1  & 2.4827 $\times 10^5$ & applied for $P_U$, 2-norm based\\
\hline
 3.4, 3.4.2  & 2.4827 $\times 10^5$ & 2-norm based\\
\hline
3.4.4, 3.4.6 & 2.4590 $\times 10^5 $ & 2-norm based\\
\hline
  3.6 & 2.4827 $\times 10^5$ & applied for $P_U$, 2-norm \\
 \hline
\end{tabular}
\vspace{0.1cm}
\caption{New upper bounds}\label{upperbound-new}
\end{table}
\end{center}
The upper bound given by Maroulas and  Psarrakos equals to $1+r_2$, with 

$  r_2= \max\{0.0059804, 0.065468, 0.84200, 0.87573, 25.012, 322.89, 322.74, 3.0513\times 10^4, 2.7181 \times 10^5\}=  2.7181 \times 10^5 .$

We can compute the maximal modulus of the eigenvalues of $P(z)$, which  is exactly $ \bold{2.4354\times 10^5}$. Moreover, Corollary \ref{coro3} and Corollary \ref{coro4} give usually the  best upper bounds.

The lower bounds obtained by Higham and Tisseur \cite{HT}  are given in Table \ref{HTlower}, while our new lower bounds are given in Table \ref{newlower}. 
\begin{center}
\begin{table}[h!]
\begin{tabular}{|l|l|l|}
\hline 
\textbf{Lemmas} & \textbf{Values} & \textbf{Comments \quad \quad \quad \quad \quad \quad }\\
\hline 
2.2&9.1316 $\times 10^{-10}$ & 2-norm\\
\hline
2.3 (2.1), 2.4 (2.5) & 8.0663 $\times 10^{-10}$ & 1-norm \\
\hline
2.3 (2.2), 2.4 (2.6) & 6.0215 $\times 10^{-10}$ & $\infty$-norm \\
\hline
2.3 (2.3), 2.4 (2.7) & 1.0456 $\times 10^{-9}$ & 2-norm \\
\hline
2.6 & 4.2286 $\times 10^{-8}$ & applied for $C_L(\alpha)$, $\beta=1/4$\\
\hline
\end{tabular}
\vspace{0.1cm}
\caption{Higham and Tisseur's lower bounds}\label{HTlower}
\end{table}
\end{center}

\begin{center}
\begin{table}[h!]
\begin{tabular}{|l|l|l|}
\hline 
\textbf{Theorems} & \textbf{Values} & \textbf{Comments \quad \quad \quad \quad \quad \quad \quad \quad \quad \quad }\\
\hline 
 3.4.1, 3.4.3  &   3.53 $\times 10^{-5}$ & 2-norm based\\
\hline
3.4.5 & 0.71005 $\times 10^{-5} $ & 2-norm based\\
\hline
 3.4.7 & 0.71034 $\times 10^{-5}$ & 2-norm based\\
 \hline
 3.5 & 9.4306 $\times 10^{-55}$ & 2-norm based\\
 \hline
 3.6 & 2.4502 $\times 10^{-10}$ & 2-norm based\\
 \hline
\end{tabular}
\vspace{0.1cm}
\caption{New lower bounds}\label{newlower}
\end{table}
\end{center}
The lower bound given by Maroulas and  Psarrakos is  $r_1 =    1.1436 \times 10^{-7}$.

We can compute the minimum modulus of the eigenvalues of $P(z)$, which is  exactly    \textbf{0.012037}. Hence the lower bounds obtained above are in general far away the expected one.  However, compare together, Corollary \ref{coro1} and Corollary \ref{coro2'} give usually the  best lower bounds.

\end{example}

In the next example we compute and compare the obtained bounds for eigenvalues of monic matrix polynomials whose coefficients are symmetric random matrices.

\begin{example} \rm 

Consider a $5\times 5$ monic matrix polynomial $P(z)$ of degree $m=9$ whose coefficient matrices are given by
$$A_i=(B_i+B_i^*)/2, ~ i=0,\ldots, 8, $$
where $B_0=B_1=rand(5)$, and $B_j=j*rand(5)$ for $j=2,\ldots,8$. Here  $rand(5)$ denotes a $5 \times 5$ random matrix from the normal $(0,1)$ distribution.

The upper bounds obtained by Higham and Tisseur \cite{HT}  are given in Table \ref{higham-upperbound}, while our new upper bounds are given in Table \ref{upperbound-new}.  
\begin{center}
\begin{table}[h!]
\begin{tabular}{|l|l|l|} 
\hline 
\textbf{Lemmas \quad } & \textbf{Values} & \textbf{Comments \quad \quad \quad \quad \quad \quad \quad \quad \quad\quad   }\\
\hline 
2.3 ~(2.2) & 99.050 & $\infty$-norm based\\
\hline
2.3 ~ (2.3) & 35.335 &2-norm based \\
\hline
2.5~ (2.13)& 31.030 & $\infty$-norm based\\
\hline
2.11~ (2.18) &  28.504  & 2-norm based\\
\hline
3.1 & 20.62502   & Cauchy's theorem applied for $P$, 2-norm \\
\hline
3.1 & 20.62502  & Cauchy's theorem applied for $P_U$, 2-norm \\
\hline
4.1& 39.542& 2-norm based\\
\hline
\end{tabular}
\vspace{0.1cm}
\caption{Higham and Tisseur's upper bounds}  \label{higham-upperbound}
\end{table}

\begin{table}[h!]
\begin{tabular}{|l|l|l|}
\hline 
\textbf{Theorems/Corollaries} & \textbf{Values} & \textbf{Comments \quad \quad \quad \quad \quad \quad \quad  \quad }\\
\hline 
3.2  & 20.77125  & 2-norm based\\
\hline
 3.2.1  & 19.77125  & 2-norm based\\
 \hline
 3.3  & 20.771  & 2-norm based\\
  \hline
 3.3.1  & 19.280  & 2-norm based\\
 \hline
 3.4 & 20.632 & 2-norm based\\
\hline
 3.4.2  & 19.455 & 2-norm based\\
\hline
3.4.4 & 19.892 & 2-norm based\\
\hline
3.4.6 & 20.697 & 2-norm based\\
\hline
  3.6 & 20.771 & applied for $P_U$, 2-norm \\
 \hline
\end{tabular}
\vspace{0.1cm}
\caption{New upper bounds}\label{upperbound-new}
\end{table}
\end{center}
The upper bound given by Maroulas and  Psarrakos equals to $1+r_2$, with 

$  r_2= \max\{1.9326, 1.3831, 4.3089, 6.3952,  5.5836,  8.3366,  19.922,  7.6824, 6.9319\}=  19.922.$

We can compute the maximal modulus of the eigenvalues of $P(z)$, which is exactly $ \bold{19.009}$. Moreover, Corollary \ref{Cauchytype4}  gives usually the  best upper bounds.

The lower bounds obtained by Higham and Tisseur \cite{HT}  are given in Table \ref{HTlower}, while our new lower bounds are given in Table \ref{newlower}. 
\begin{center}
\begin{table}[h!]
\begin{tabular}{|l|l|l|}
\hline 
\textbf{Lemmas} & \textbf{Values}\quad \quad  & \textbf{Comments \quad \quad \quad \quad \quad \quad \quad \quad \quad  }\\
\hline 
2.2&0.0027435 & 2-norm\\
\hline
2.3 (2.1), 2.4 (2.5) & 0.0084515 & 1-norm \\
\hline
2.3 (2.2), 2.4 (2.6) & 0.0019335 & $\infty$-norm \\
\hline
2.3 (2.3), 2.4 (2.7) & 0.0067592 & 2-norm \\
\hline
2.6 & 0.0093639 & applied for $C_L(\alpha)$, $\beta=1/4$\\
\hline
\end{tabular}
\vspace{0.1cm}
\caption{Higham and Tisseur's lower bounds}\label{HTlower}
\end{table}
\end{center}

\begin{center}
\begin{table}[h!]
\begin{tabular}{|l|l|l|}
\hline 
\textbf{Theorems} & \textbf{Values} & \textbf{Comments \quad \quad \quad \quad \quad \quad \quad \quad \quad \quad }\\
\hline 
 3.4.1  &   0.073516 & 2-norm based\\
\hline
3.4.3 & 0.068055  & 2-norm based\\
\hline
 3.4.5 & 0.058889 & 2-norm based\\
 \hline
 3.5 & 3.9543  $\times 10^{-15}$ & 2-norm based\\
 \hline
 3.6 & 0.0024740 & 2-norm based\\
 \hline
\end{tabular}
\vspace{0.1cm}
\caption{New lower bounds}\label{newlower}
\end{table}
\end{center}
The lower bound given by Maroulas and  Psarrakos is  $r_1 =    1.3606$.

We can compute the minimum modulus of the eigenvalues of $P(z)$, which is  exactly      \textbf{0.15023}.  Compare together, Corollary \ref{coro1} and Corollary \ref{coro2'} give usually the  best lower bounds.
\end{example}
\section*{Acknowledgements}
The authors would like to thank the anonymous referees for their useful comments and suggestions to improve the results in the original version of this paper.  The final version of this paper was finished during the visit of the first author   at the Vietnam Institute for Advanced Study in Mathematics (VIASM). He thanks VIASM for financial support and hospitality.

\end{document}